\documentclass[final,twoside,web]{ieeecolor}

\def\BibTeX{{\rm B\kern-.05em{\sc i\kern-.025em b}\kern-.08em
   T\kern-.1667em\lower.7ex\hbox{E}\kern-.125emX}}
\usepackage{flushend}

\usepackage{generic}
\usepackage{cite}
\usepackage[noend]{algorithmic}
\usepackage[ieee]{mathematix}

\usepackage{filecontents}

\usepackage[T1]{fontenc}
\usepackage[utf8]{inputenc}

\usepackage{booktabs}
\usepackage{multirow}
\usepackage{multicol}
\usepackage{siunitx} 

\usepackage{amssymb,amsmath,mathtools}
\usepackage{units}
\usepackage{xspace}
\usepackage{csquotes}
\usepackage{textcomp}
\usepackage{cleveref}
\usepackage{url}

\usepackage{tikz}
\usepackage{pgfplots}
\pgfplotsset{compat=1.9}
\usepgfplotslibrary{statistics}

\usepgfplotslibrary{external}
\usepgfplotslibrary{groupplots}
\tikzset{external/system call={pdflatex \tikzexternalcheckshellescape -halt-on-error
			-interaction=batchmode -jobname "\image" "\texsource"}}
\usetikzlibrary{matrix}
\usetikzlibrary{patterns}
\usetikzlibrary{arrows}
\usetikzlibrary{decorations.pathreplacing,decorations.markings}
\usetikzlibrary{backgrounds}
\usetikzlibrary{calc}

\usepackage{pgfplotstable}
\usetikzlibrary{patterns}

\tikzset{ 
	on each segment/.style={
			decorate,
			decoration={
					show path construction,
					moveto code={},
					lineto code={
							\path [#1]
							(\tikzinputsegmentfirst) -- (\tikzinputsegmentlast);
						},
					curveto code={
							\path [#1] (\tikzinputsegmentfirst)
							.. controls
							(\tikzinputsegmentsupporta) and (\tikzinputsegmentsupportb)
							..
							(\tikzinputsegmentlast);
						},
					closepath code={
							\path [#1]
							(\tikzinputsegmentfirst) -- (\tikzinputsegmentlast);
						}
				}
		}
}

\tikzset{ 
	mid arrow/.style={
			postaction={
					decorate,decoration={
							markings,
							mark=at position .5 with {\arrow[#1]{stealth}}
						}
				}
		}
}

\newcommand{\minfbe}{\textsc{minfbe}\xspace}
\newcommand{\nama}{\textsc{nama}\xspace}
\newcommand{\gpad}{\textsc{gpad}\xspace}

\newcommand{\smallplus}{{\scriptscriptstyle +}}

\usepackage{color}
\definecolor{plecs1}{rgb}{0,.7969,0}
\definecolor{plecs2}{rgb}{1,0,0}
\definecolor{plecs3}{rgb}{0,0,1}
\definecolor{plecs4}{rgb}{1,.7969,0}
\definecolor{plecs5}{rgb}{1,0,1}
\definecolor{plecs6}{rgb}{0,1,1}

\definecolor{vgRed}{RGB}{193, 48, 24}
\definecolor{vgOrange}{RGB}{243, 111, 19}
\definecolor{vgYellow}{RGB}{235, 203, 56}
\definecolor{vgGreen}{RGB}{162, 185, 105}
\definecolor{vgLightBlue}{RGB}{13, 149, 188}
\definecolor{vgDarkBlue}{RGB}{6, 56, 81}

\usepackage[acronym]{glossaries}

\usepackage{acronym}
\acrodef{gpu}[GPU]{Graphic processing unit}
\acrodef{mpc}[MPC]{model predictive control}
\acrodef{admm}[ADMM]{alternating direction of method multipliers}

\begin{document}
\title{Massively parallelizable proximal algorithms for large-scale stochastic optimal control problems}
\author{Ajay~K.~Sampathirao,
   Panagiotis~Patrinos, \IEEEmembership{Member, IEEE},
   Alberto~Bemporad, \IEEEmembership{Fellow, IEEE},
   and Pantelis~Sopasakis
   \thanks{A.K. Sampathirao is with Enervalis NV, Belgium,
      formerly with Technical University of Berlin (TUB),
      Control Systems Lab (e-mail: \texttt{Sampathirao@control.TU-Berlin.de}).}
   \thanks{P. Patrinos is with KU Leuven, Department of Electrical Engineering (ESAT),
      (e-mail: \texttt{panos.patrinos@esat.kuleuven.be}).}
   \thanks{A. Bemporad is with IMT School for Advanced Studies Lucca, Italy,
      (email: \texttt{alberto.bemporad@imtlucca.it}).}
   \thanks{P. Sopasakis is with Queen's University Belfast, School of Electronics, Electrical
      Engineering and Computer Science (EEECS) and the Centre for Intelligent Autonomous Manufacturing
      Systems (\textit{i}-AMS), Ashby Building, Stranmillis Road, BT9 5AG,
      United Kingdom (e-mail: \texttt{p.sopasakis@qub.ac.uk}).}}

\maketitle

\begin{abstract}
   Scenario-based stochastic optimal control problems suffer from the curse of dimensionality
   as they can easily grow to six and seven figure sizes. First-order methods are suitable
   as they can deal with such large-scale problems, but may fail to achieve accurate solutions
   within a reasonable number of iterations. To achieve solutions of higher accuracy and high
   speed, in this paper we propose two proximal quasi-Newtonian limited-memory
   algorithms
   ---
   \minfbe applied to the dual problem and
   the Newton-type alternating minimization algorithm (\nama)
   ---
   which can be massively parallelized on lockstep hardware such as
   graphics processing units (GPUs).
   We demonstrate the performance of these methods, in terms of
   convergence speed and parallelizability, on large-scale problems
   involving millions of variables.
\end{abstract}


\section{Introduction}
\label{sec:introduction}

\subsection{Background}
\IEEEPARstart{S}{tochastic} optimal control is the backbone of
stochastic \ac{mpc}, which is known for its appealing stability and
constraint satisfaction properties~\cite{PatPanHar+2014, ChaLyg2015}
and has found several applications~\cite{LiaSixCha2016, CaiBerBem+2014,
  DarGeoSmiLyg2019}.
More specifically, scenario-based stochastic \ac{mpc} is gaining great
popularity~\cite{HanSopBemRaiCol15,SamSopBem+2018,ZidKolBem2021,MatVirKol+2021}
due to its applicability to virtually any stochastic model of uncertainty that can
be reasonably approximated by a discrete distribution.
However, the limiting factor towards its industrial uptake is the
computational time required to solve numerically the resulting large-scale
optimisation problem.
Indeed, multistage scenario-based stochastic optimal control problems suffer from
the curse of dimensionality and can lead to problems with millions of
decision variables \cite{SamSopBem+2018}.

\acp{gpu} have been used for their massive parallelization capabilities
in applications as diverse as
cryptocurrency mining \cite{bitcoin-mining},
cosmology \cite{aubert2010numerical},
medical image processing \cite{eklund2013medical}, simulations of molecular
dynamics \cite{le2013spfp}, machine learning \cite{kim2017performance},
and a lot more.  \acp{gpu} are suitable for lockstep parallelization,
where the same elementary operations are applied to different memory
positions using dedicated functions known as \textit{kernels}.
Programming \acp{gpu} for general-purpose data-parallel computations
is facilitated by programming languages and frameworks such as
CUDA \cite{cuda, CookCuda2012} (for NVIDIA GPUs, used by well-known
software such as Tensorflow \cite{tensorflow2015-whitepaper} and Caffe
\cite{jia2014caffe}), OpenCL, OpenACC, OpenGL and more.

In recent years, a number of papers have proposed parallelizable
variants of numerical optimization methods such as the interior
point method~\cite{GadeNielsen14}, parallel quadratic
programming~\cite{Yu2017}, \ac{admm}~\cite{Fang18,Enfedaque+18,
  QureshiEastCannon2019}
and other proximal algorithms~\cite{GaetanoChierchia2012,Schubiger19}.
In these approaches, \acp{gpu} are used to parallelize the involved
algebraic operations and the solution of linear systems:
the primal-dual optimalily conditions in interior point algorithms
and equality-constrained QPs in \ac{admm}.
Given the lockstep data parallelization paradigm of \acp{gpu},
numerical methods that aim at splitting the problem into smaller
optimization problems that are to be executed in parallel
(such as \cite{Deng2016} and \cite{kouzoupis2019dual}) do not lend
themselves to \ac{gpu} implementations.

Scenario-based problems possess a certain structure that
can be exploited to design very efficient \textit{ad hoc} GPU-enabled
implementations leading to a higher acceleration as discussed in
\cite{SamSopBem+2018}.
It has been shown that first-order algorithms such as the accelerated
proximal gradient method can be used to achieve significant speed-ups
\cite{sopasakis2018uncertainty,SamSopBem+2018,SamSopBemPat15}.
However, first-order methods tend to be prone to ill-conditioning as
they disregard curvature information.
This motivates the development of numerical methods that can exploit
the underlying problem structure of scenario-based optimal control problems,
come with good convergence characteristics, and are amenable to lockstep
parallelisation on \acp{gpu}.

In this paper we propose two massively parallelizable numerical methods
that exploit the structure of scenario-based stochastic optimal control
problems, building up on
(i) the \minfbe method \cite{TheStePat16} applied to the dual problem,
(ii) the Newton-type alternating minimization algorithm (\nama) \cite{nama2019} algorithms,
as well as (iii) on our previous work on GPU-accelerated optimization \cite{SamSopBemPat17a}.
All methods lend themselves to highly
parallelizable implementations and lead to similar convergence speeds.
However, we will show that \nama allows a significantly higher parallelizability
and lower computation times.
\minfbe and \nama involve only simple algebraic operations, use limited-memory BFGS
directions and can achieve better accuracy and significantly faster
convergence than the accelerated proximal gradient method of \cite{SamSopBem+2018}
(linear convergence rate instead of $\mathcal{O}(1/k^2)$).

\subsection{Notation}
Let $\N$, $\R$, $\R^n$ and $\R^{m{}\times{}n}$ denote the sets of nonnegative integers,
real numbers, $n$-dimensional vectors and $m$-by-$n$ matrices respectively.
Let $\N_{[k_1, k_2]} \dfn \{n\in\N{}:{} k_1\leq{}n{}\leq{}k_2\}$. Let $\barre=\R{}\cup{}\{+\infty\}$
be the set of extended-real numbers. Given a set $X{}\subseteq{}\R^n$ and $x\in\R^n$ we define the
\textit{indicator} of $X$ as the extended-real-valued function $\delta({}\cdot{}\mid{}X):\R^n\to\barre$
with $\delta(x{}\mid{}X) = 0$ for $x{}\in{}X$ and $\delta(x {}\mid{} X)=\infty$ otherwise.
For $A\in\R^{m{}\times{}n}$, $A^\top$ denotes the transpose of $A$.
For $A,B\in\R^{m\times n}$, we write $A\succ{}B$ ($A \succcurlyeq B$) if $A-B$ is positive
(semi)definite. For a convex function $f:\R^n\to\barre$, its convex conjugate function $f^*$
is defined as
\(
f^*(y) = \sup_{x}\{x^\top y-f(x)\}.
\)
Lastly, given a nonempty, closed, convex set $X{}\subseteq{}\R^n$, we define the projection
operator onto $X$ as $\proj_X(x) = \argmin_{y{}\in{}X}\|y-x\|$.


\section{Problem Statement}
We start by stating the stochastic optimal control problem
we will study in this paper.

\subsection{Stochastic dynamics on scenario trees}
\begin{figure}[h]
    \centering
    \includegraphics[width=0.98\linewidth]{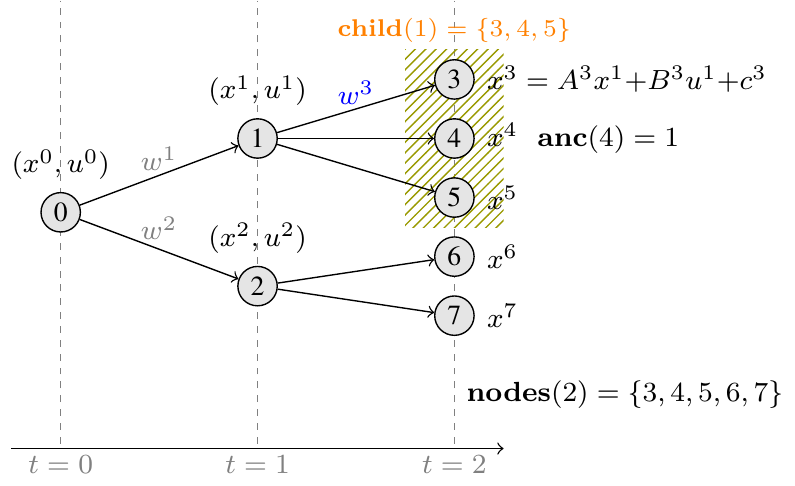}
    \caption{Scenario tree structure with three stages and the system dynamics
        on its nodes. The cost associated with node $i=1$ is $\ell_1(x^0, u^0, w^1)$.}%
    \label{fig:scenario_tree}
\end{figure}

Consider a discrete-time stochastic dynamical system
of the form
\begin{equation}\label{eq:affine_stochastic_dynamics}
    x_{t+1} = A_{w_t}x_t + B_{w_t}u_t + c_{w_t},
\end{equation}
with state $x_t\in\R^{n_x}$ and input $u_t\in\R^{n_u}$, which is
driven by the stochastic process $w_t$.
For example, Markov jump affine systems fall into this category \cite{bemporad+2018}.
The evolution of this system over a finite sequence
of time instants, $t\in\N_{[0, N]}$,
can be described using a \textit{scenario tree}: a directed graph of the
form shown in \Cref{fig:scenario_tree}.
The scenario tree structure is essentially the representation of a
discrete multistage probability distribution.
A scenario tree represents the evolution of the system states as
more information becomes available: at every stage $t$, we
assume that the state, $x_t$, can be measured and a control
action $u_t$ can be decided based on that measurement, thus modeling
an entire feedback policy.

The nodes of the scenario tree are organised in \textit{stages}, $t\in\N_{[0, N]}$,
and indexed by a unique integer $i$. At stage $t=0$ we assume that the state
--- which is the current state in an \ac{mpc} setting ---
is known; this corresponds to the \textit{root} node of the tree, which is
indexed by $i=0$. The nodes at a stage $t$ are denoted by $\nodes(t)$
and the nodes at stage $t=N$ are called the \textit{leaf nodes} of the tree.
For notational convenience, we will denote the nodes at stages $t\in\N_{[t_1, t_2]}$,
with $0 \leq t_1 \leq t_2 \leq N$,
by $\nodes(t_1, t_2) = \bigcup_{t=t_1}^{t_2}\nodes(t)$.
The set $\nodes(t)$ is a probability space: every node $i\in\nodes(t)$
is assigned a nonzero probability value $\pi^i$. Naturally $\pi^0 = 1$
and $\sum_{i\in\nodes(t)}\pi^i = 1$ for all $t\in\N_{[0, N]}$.

Every node $i$ at a stage $t\in\N_{[1, N]}$ has an \textit{ancestor},
$\anc(i) \in \nodes(t-1)$, and all nodes at a stage $t\in\N_{[0, N-1]}$
have a set of \textit{children}, $\child(i) \subseteq \nodes(t+1)$.
The set $\child(i)$ is a probability space with probability vector
$\pi^{[i]} \in \R^{|\child(i)|}$. This is a vector whose $i_\smallplus$-th
element is equal to ${\pi^{i_\smallplus}}/{\pi^i}$ --- for short
\(
\pi^{[i]}
{}={}
\tfrac{1}{\pi^i}(\pi^{i_{\smallplus}})_{i_{\smallplus}\in\child(i)}.
\)

The system dynamics, \eqref{eq:affine_stochastic_dynamics}, across the
nodes of the scenario tree can be stated as
\begin{equation}\label{eq:system_dynamics_tree}
    x^{i_\smallplus}
    {}={}
    A^{i_\smallplus}x^{i}
    {}+{}
    B^{i_\smallplus} u^{i}
        {}+{}
    c^{i_\smallplus},
\end{equation}
for $i\in\nodes(0, N-1)$, $i_\smallplus\in\child(i)$.
Note that the total number of scenarios coincides with the number of
leaf nodes, and the number of non-leaf nodes with the number of free
input variables (see \Cref{fig:scenario_tree}).


\subsection{Stochastic optimal control problem}
A \textit{multistage stochastic optimal control problem} for~\eqref{eq:affine_stochastic_dynamics}
with horizon $N$ can be formulated as
%
%
%
\begin{align*}
    \mathbb{P}(p)
    {}:{}
    \minimize_{\{u_t\}_{t=0}^{N-1},\{x_t\}_{t=0}^{N}}\
    \mathbb{E}\left[V_f(x_N) + \sum_{t=0}^{N-1}\ell_t(x_t, u_t,w_t)\right],
\end{align*}
subject to~\eqref{eq:affine_stochastic_dynamics} and the condition
$x_0=p$.
Note that in this formulation, $\{u_t\}_{t=0}^{N-1}$ and $\{x_t\}_{t=0}^{N}$
are random variables.
The \textit{stage cost} at stage $t\in\N_{[1, N]}$ is a random
variable which admits the values $\ell^{i}(x^{\anc(i)}, u^{\anc(i)}) \dfn \ell_t(x^{\anc(i)}, u^{\anc(i)}, w^{i})$,
for $i\in\nodes(t)$, with probability $\pi^{i}$.
The \textit{terminal cost} function is also a random variable which
admits the values $V_f(x^i)$ for $i\in\nodes(N)$ with probability $\pi^i$.
That said, the optimal control problem can be written as
\begin{multline*}
    \mathbb{P}(p){}:{} \minimize_{
    \substack{
    \{u^i\}_{i\in\nodes(0, N-1)}
    \\
    \{x^i\}_{i\in\nodes(0, N)}
    }
    \hspace{0.7em}
    }
    \sum_{
        i\in\nodes(1, N)
    }
    \hspace{-1.5em}
    \pi^i\ell^i(x^{\anc(i)}, u^{\anc(i)})
    \\
    +\sum_{i\in\nodes(N)}\hspace{-1em}\pi^i V_f^i(x^i),
\end{multline*}
subject to the system dynamics \eqref{eq:system_dynamics_tree} and the condition
$x_0=p$.

The stage cost function, $\ell^{i}:\R^{n_x}\times\R^{n_u}\to\barre$,
at node $i\in\nodes(1, N)$, is an extended-real-valued
function which can be decomposed as follows
\begin{equation}
    \ell^i(x, u) = \phi^i(x, u) + \bar{\phi}^i(F^i x + G^i u),
\end{equation}
where $\phi^i:\R^{n_x}\times\R^{n_u}\to\R$ is a smooth convex function and
$\bar{\phi}^i:\R^{m_i}\to\barre$ is a proper, extended-real-valued, possibly
nonsmooth, convex, lower semicontinuous function and
$F^i \in \R^{m_i\times n_x}$, $G^i\in\R^{m_i\times n_u}$.
Functions $\bar{\phi}^i$ can be taken to be indicator functions so as to model
constraints on inputs and states.

We can also decompose the terminal cost function, $V_f^i:\R^{n_x}\to\barre$,
as follows
\begin{equation}
    V_f^i(x) = \phi_N^i(x) + \bar{\phi}_N^i(F_N^i x),
\end{equation}
where $F_N^i{}\in{}\R^{m_{N,i}{}\times{}n}$,
$\phi_N^i:\R^{n_x}\to\R$ is real valued, smooth, convex function and
$\bar{\phi}_N^i:\R^{m_{N,i}}\to\barre$ is a proper extended-real-valued, convex,
lower semicontinuous function.

Functions $\bar{\phi}^{i}$ need not be smooth. They can be used to
describe hard joint state-input constraints of the form
$F^{i} x + G^{i} u \in Y^{i}$ by taking $\bar\phi^i(\cdot)=\delta({}\cdot{}|{}Y^i)$.
Similarly, $\bar{\phi}^{i}$ can describe soft constraints simply
by replacing the indicator function $\delta({}\cdot{}|{}Y^i)$ by a distance-to-set
function. On the other hand, functions $\phi^i$ and $\phi_N^i$ are typically
taken to be convex quadratic (and $\phi^i$ are assumed to be strongly convex
with respect to $u$ and jointly convex in $(x, u)$).
Hereafter, we consider the quadratic cost functions
\begin{equation}
    \phi^i(x,u)= \begin{bmatrix}x\\u\end{bmatrix}^\top
    \begin{bmatrix}Q_i & S_i^\top\\ S_i & R_i\end{bmatrix}
    \begin{bmatrix}x\\u\end{bmatrix} +
    q_i^\top x +r_i^\top u,
\end{equation}
for $x\in\R^{n_x}$ and $u\in\R^{n_u}$,
with $Q_i = Q_i^\top \succcurlyeq 0$, and $R_i=R_i^\top \succ 0$,
and
\begin{equation}\label{eq:95395e1a-21ec-4e2c-b835-d6ab88f12dc8}
    \begin{bmatrix}
        Q_i & S_i^\top \\ S_i & R_i
    \end{bmatrix}
    {}\succcurlyeq{}
    0,
\end{equation}
for all $i\in\nodes(0, N-1)$.
Lastly, $V_f(x) = x^\top P_N x + p_N^\top x$ for $x\in\R^{n_x}$ with $P_N = P_N^\top \succ 0$.

\subsection{Formulation of optimization problem}
The decision variable of $\mathbb{P}(p)$ is the vector
$x=\left((u^i)_{i\in\nodes(0, N-1)}, (x^{i})_{i\in\nodes(1, N)}\right)\in\R^n$,
where $n=|\nodes(0, N-1)|n_u + |\nodes(1, N)|n_x$.
Let us define the affine space
\begin{align}
    \mathcal{Z}(p) = \left\{x{}\left|
    \begin{array}{l}
        x^0 = p, x^{i_\smallplus} = A^{i^\smallplus} x^i + B^{i_\smallplus}u^{i} + c^{i_\smallplus},
        \\
        i\in\nodes(0, N-1), i_\smallplus \in \child(i).
    \end{array}
    \right.\right\},
\end{align}
which describes the system dynamics.
Let us also define the functions $f(\cdot;p):\R^n\to\barre$
and $g{}:{}\R^m\to\barre$ that maps $z=((z^i)_{i\in\nodes(0,N)},
    (z_{N}^{i})_{i\in\nodes(N)})$ with $z^i\in\R^{m_i}$ and
$z_{N}^{i}\in\R^{m_{N,i}}$ and is given by
%
%
\begin{subequations}
    \begin{align}
        f(x) & = \sum_{i\in\nodes(1, N)}\hspace{-1em}
        \pi^i\phi^i(x^{\anc(i)},u^{\anc(i)})\notag
        \\
             & \qquad {}+{}\hspace{-1em}
        \sum_{i\in\nodes(N)}\hspace{-1em}\pi^i \phi_N(x^i)
        {}+{}
        \delta(x| \mathcal{Z}(p)),
        \label{eq:f-definition}
        \\
        g(z) & {}={}\hspace{-0.8em}
        \sum_{i\in\nodes(1, N)}\hspace{-1em}
        \pi^i\bar\phi^i(z^i)
        {}+{}\hspace{-0.4em}
        \sum_{i\in\nodes(N)}\hspace{-1em}
        \hspace{-0.4em}
        \pi^i\bar\phi_N(z_N^i),
    \end{align}
\end{subequations}
where $z=((z^i)_{i\in\nodes(0, N)}, (z_N^i)_{i\in\nodes(N)})$ and
define $H:\R^n\to\R^m$ as a linear operator that maps $x$ to a vector $z\in\R^m$
as above with $z^i = F^i x^{\anc(i)} + G^i u^{\anc(i)}$ for $i\in\nodes(0, N)$
and $z_N^i = F_N^i x^i$ for $i\in\nodes(N)$.

Given that functions $\phi^{i}$ are quadratic as described in the previous section,
function $f$ is strongly convex (as it follows from \cite[Prop. 6]{frison2012msc}), therefore the convex conjugate of $f$, $f^*$,
is differentiable with $L$-Lipschitz gradient because of~\cite[Prop.~{12.60}]{rockafellar-book}.

Problem $\mathbb{P}(p)$ can be written as
\begin{align}\label{eq:p-star}
    \mathbb{P}(p){}:{} \minimize_{x\in\R^n} f(x; p) + g(Hx).
\end{align}
Hereafter, we assume that $\mathbb{P}(p)$ is feasible.
The Fenchel dual of Problem $\mathbb{P}(p)$  in \Cref{eq:p-star} is
\begin{align}\label{eq:dual-problem}
    \mathbb{D}(p){}:{} \minimize_{y\in\R^m} f^*(-H^\top y; p) + g^*(y).
\end{align}
Let us define the function $\hat{f}:\R^m\to\R$ as
\begin{equation}
    \hat{f}(y; p)\dfn f^*(-H^\top y; p).
\end{equation}
Then, Problem $\mathbb{D}(p)$ can be written as
\begin{equation}
    \mathbb{D}(p){}:{} \minimize_{y\in\R^m} \hat{f}(y; p) + g^*(y).
\end{equation}
For given $p$, strong duality holds if there is an
$x\in\mathcal{Z}(p)$ such that
$Hx \in \operatorname{relint} \dom g$ \cite[Theorem 15.23]{BauschkeCombettes2011}
--- we will hereafter assume that this assumption is satisfied.


\subsection{Optimality conditions}

The proximal operator of a proper, closed, convex function $g$
plays a major role in modern optimization theory and is defined as
\begin{align}
    \prox_{\lambda g}(v) = \argmin_{z}\{g(z) + \tfrac{1}{2\lambda}\|v-z\|^2\},
\end{align}
with $\lambda{}>{}0$.
Proximal operators of a great variety of functions including indicators
of sets, distance-to-set functions and norms can be easily evaluated
analytically and at a very low computational cost~\cite{ComPes10}.
For example, the proximal operator of the indicator of a set $Y$
is the projection on $Y$, that is
$\prox_{\lambda \delta(\cdot\mid Y)}(v)=\proj(v\mid Y)$.

A simple optimality condition for~\eqref{eq:dual-problem} is
\begin{align}\label{eq:optimality-condition}
    y-\prox_{\lambda g^*}(y-\lambda \nabla \hat{f}(y)) = 0,
\end{align}
for some $\lambda>0$~\cite{parikh-boyd-proximal}.
By virtue of the Moreau decomposition formula,~\eqref{eq:optimality-condition}
is equivalently written as
\begin{align}
    \nabla \hat{f}(y) + \prox_{\lambda^{-1}g}(\lambda^{-1} y - \nabla \hat{f}(y)) = 0.
\end{align}

We define the \textit{forward-backward mapping}
\begin{align}
    T_\lambda(y) & \dfn \prox_{\lambda g^*}(y-\lambda \nabla \hat{f}(y)),
\end{align}
which, using the Moreau decomposition property, becomes
\begin{align}\label{eq:t_lambda_prox}
    T_\lambda(y) = y - \lambda \nabla \hat{f}(y) - \lambda \prox_{\lambda^{-1}g}(\lambda^{-1}y {-} \nabla \hat{f}(y)),
\end{align}
and we also define the \textit{fixed-point residual} mapping
\begin{align}
    R_\lambda(y) & \dfn \lambda^{-1}(y-T_\lambda(y))      \\
                 & = z_{\lambda}(y)-Hx(y),\label{eq:fpr2}
\end{align}
where $x(y)$ and $z_\lambda(y)$ are defined as
\begin{subequations}
    \begin{align}
        x(y)           & \dfn \nabla f^*(-H^\top y),                                 \\
        z_{\lambda}(y) & \dfn \prox_{\lambda^{-1}g}\left(\lambda^{-1}y+Hx(y)\right),
    \end{align}
\end{subequations}
therefore,
\begin{equation}
    x(y)=\argmin_{z}\{\< z,H^\top y\> + f(z)\}.
\end{equation}
Note also that $T_\lambda(y)$ can be computed from \Cref{{eq:t_lambda_prox}} as
\begin{equation}
    T_\lambda(y)
    =
    y - \lambda \left(\nabla \hat{f}(y) - z_\lambda(y)\right).
\end{equation}

The aforementioned optimality condition in~\Cref{eq:optimality-condition} is
equivalently written as $R_\lambda(y)=0$, that is, solving the dual optimization
problem~\eqref{eq:dual-problem} becomes equivalent to finding a zero of the
operator $R_\lambda$.


\section{Numerical Optimization}

\subsection{Forward-backward envelope}
The \textit{forward-backward envelope} (FBE) of~\eqref{eq:dual-problem}
is a real-valued function $\varphi_\lambda$ given by~\cite{patrinos-FB-TruncNewton-2014,patrinos-bemporad-2013}
\begin{align}
  \varphi_\lambda(y) = & \hat{f}(y)
  + g^*(T_\lambda(y))
  \notag                            \\&- \lambda \<  \nabla \hat{f}( y), R_\lambda(y)\>
  + \tfrac{\lambda}{2}\|R_\lambda(y)\|^2.
\end{align}
If $\hat{f}$ is twice continuously differentiable --- conditions under which this is the case
can be found in \cite{Gorni1991} --- then $\varphi_\lambda$ is continuously
differentiable with
\begin{equation}
  \nabla \varphi_\lambda(y) = (I-\lambda \nabla^2 \hat{f}(y))R_\lambda(y).
\end{equation}
Note that in practice it is not necessary to compute or store the Hessian matrix
$\nabla^2 \hat{f}(y)$. Instead, it suffices to implement an algorithm that returns
Hessian-vector products of the form $\nabla^2 \hat{f}(y){}\cdot{}z$.
The most important property of the FBE is that for $\lambda\in(0,\nicefrac{1}{L})$,
the set of minimizers of~\eqref{eq:dual-problem} coincides with
\begin{align*}
  \argmin \varphi_\lambda & \equiv \zer \nabla \varphi_{\lambda}\dfn \{y: \nabla \varphi_\lambda(y) = 0 \} \\
                          & = \argmin \hat{f}(y) + g^*(y)
  =\zer R_\lambda.
\end{align*}

Essentially, the problem of solving the dual optimization problem~\eqref{eq:dual-problem}
is equivalent to the unconstrained minimization of the continuously differentiable function $\varphi_\lambda$,
that is
\begin{subequations}\label{eq:equivalence_fbe_minimization}
  \begin{align}
    \inf \hat{f}(y)  + g^*(y) {}={}   & \inf \varphi_\lambda,    \\
    \argmin \hat{f}(y) + g^*(y) {}={} & \argmin \varphi_\lambda.
  \end{align}
\end{subequations}
Moreover, the above is equivalent to finding a zero of the fixed-point residual operator.
In the common case where $\hat{f}$ is strongly convex quadratic, $\phi_\lambda$ is both
continuously differentiable and convex.

\subsection{Dual \minfbe method}
If $\hat{f}$ is twice differentiable, according to \eqref{eq:equivalence_fbe_minimization}
the original (dual) optimization problem can be cast as an \textit{unconstrained} optimization
problem with a smooth cost function. As a result we can use an appropriate unconstrained
optimization method to solve such problems, such as limited-memory BFGS \cite{SamSopBemPat17a},
however, convergence is only guaranteed under restrictive requirements
(such as twice differentiability and uniform convexity of the FBE~\cite{SunYuan2006}).

Instead, \minfbe is a method that can be applied to problems with nonsmooth
cost functions using the forward-backward envelope as a merit function using
a simple line search~\cite{TheStePat16}. \minfbe involves simple and computationally
inexpensive iterations, and exhibits superior global convergence properties.
The application of \minfbe to the dual optimization problem, $\mathbb{D}(p)$,
leads to \Cref{alg:FB-LBFGS}.
\begin{algorithm}[H]
  \caption{Dual \minfbe with L-BFGS directions}\label{alg:FB-LBFGS}
  \begin{algorithmic}[1]
    \REQUIRE $\lambda\in(0,\nicefrac{1}{L})$,  $y^0$, $m$ (memory), $\epsilon$ (tolerance)
    \ENSURE Primal-dual solution triple $(x, z, y)$
    \STATE Initialize an L-BFGS buffer with memory $m$
    \WHILE{$\|R_{\lambda}(y^k)\|_{\infty} {}>{} \epsilon$}
    \STATE {$d^k = -B^k \nabla \varphi_{\lambda}(y^k)$} (Compute an L-BFGS direction using the L-BFGS buffer)
    \STATE Choose the smallest $\tau_k \in \{2^{-\nu}\}_{\nu\in\N}$ so that
    \begin{equation}\label{eq:decrease_condition_fbe}
      \varphi_{\lambda}(w^{k})\leq \varphi_{\lambda}(y^k),
    \end{equation}
    where $w^{k} = y^k + \tau_k d^k$
    \STATE {$x^{k+1} = T_\lambda(w^k)$}
    \STATE Compute the L-BFGS-related quantities $s^k = y^{k+1} - y^k$,
    $q^k = \nabla \varphi_{\lambda}(y^{k+1}) - \nabla \varphi_{\lambda}(y^{k})$
    and $\rho_k = \langle s^k, q^k \rangle$
    \IF{$\<s^k, q^k\> > \epsilon' \|s^k\|^2 \|\nabla \varphi_\lambda(y^k)\|^2$}
    \STATE Push $(s^k,q^k,\rho^k)$ into the LBFGS buffer
    \ENDIF
    \STATE $k = k + 1$
    \ENDWHILE
    \STATE \textbf{return} $(x, z, y) = (x(y^k), z(y^k), y^k)$
  \end{algorithmic}
\end{algorithm}

\minfbe consists in applying the forward-backward mapping on the extrapolated
vector $w^k = y^k + \tau_k d^k$ which satisfies the decrease condition
\eqref{eq:decrease_condition_fbe}. The L-BFGS buffer is updated with
the vectors $s^k$, $q^k$, and their inner product $\rho_k$, provided that
the minimum-curvature condition in line 7 is satisfied for a small tolerance
$\epsilon'{}>{}0$, following \cite{Li2001}.

The algorithm iterates on the dual vectors $y^k$ and
returns a triple $(x, z, y)$ which satisfies the termination
condition $\|R_{\lambda}(y^k)\|_{\infty} \leq \epsilon$, which, in light of
\Cref{eq:fpr2} means that
\begin{subequations}
  \begin{align}
    \|z - Hx\|_{\infty} \leq \epsilon,
    \\
    -H^\top y \in \partial f(x),
    \\
    \dist_{\|\cdot\|_\infty}(y, \partial g(z)) \leq \lambda \epsilon,
  \end{align}
\end{subequations}
where \(\dist_{\|\cdot\|_\infty}\) denotes the point-to-set distance with
repsect to the $\infty$-norm.

We should highlight that the line search in line 4 of \Cref{alg:FB-LBFGS}
is a simple descent condition on the FBE, which is simpler than the Wolfe
conditions used in \cite{SamSopBemPat17a}. Moreover, although in \Cref{alg:FB-LBFGS}
we use L-BFGS directions, the method works with any direction of descent $d^k$
with respect to the FBE, that is, if $\<d^k, \nabla \varphi_\lambda(y^k)\> \leq 0$.

If $f$ is quadratic plus the indicator of an affine subspace, $x(y)$ turns
out to be linear, that is
\begin{equation}
  x(w) = x(y{}+{}\tau d) = x(y) + \tau x(d),
\end{equation}
and $\hat{f}$ is a quadratic function, that is $\nabla \hat{f}$ is linear and
\(
\hat{f}(y) = \<y, \nabla \hat{f}(y)\>,
\)
from which we can see that
\begin{align}
  \hat{f}(y+\tau d) {}={} & \<y + \tau d, \nabla \hat{f}(y+ \tau d)\>\notag
  \\
  {}={}                   & \hat{f}(y) + \tau^2 \hat{f}(d) + 2\tau\<y, \nabla\hat{f}(d)\>.
\end{align}
By virtue of the last two properties and after some algebraic manipulations,
we find that the line search condition $\varphi_{\lambda}(w^{k}){}-{}\varphi_{\lambda}(y^k) {}\leq {}0$
is equivalent to
\begin{equation}\label{eq:linesearch_trick}
  \alpha_2(y, d) \tau^2 + \alpha_1(d) \tau + \alpha_0(\tau; y, d) \leq 0,
\end{equation}
where
\begin{subequations}
  \begin{align}
    \alpha_0(\tau; y, d) {}={} & g^*(T_\lambda(y+\tau d)) - g^*(T_\lambda(y)) \notag                   \\
                               & + \tfrac{\lambda}{2}[\|z_\lambda(y+\tau d)\|^2 - \|z_\lambda(y)\|^2],
    \\
    \alpha_1(y,d) {}={}        & \<Hx(d), 2y - \lambda Hx(y)\>,
    \\
    \alpha_2(d) {}={}          & \hat{f}(d) - \tfrac{\lambda}{2}\|Hx(d)\|^2,
  \end{align}
\end{subequations}
and $\hat{f}(d)$ can be computed by invoking
\cite[Theorem 23.5]{rockafellar-convex-analysis},
from which
\begin{equation}
  \hat{f}(d) = - \<Hx(d), d\> - f(x(d)).
\end{equation}

Note that $\alpha_1$ and $\alpha_2$ do not depend on $\tau$, therefore,
can be computed once per iteration. This leads to a significant reduction of
the involved floating point operations per iteration.
The most computationally demanding parts of \minfbe are (i) the computation of
$x(y)$ and $x(d)$, and (ii) the computation of the Hessian-vector product
$\nabla\hat{f}(y^k)R_\lambda(y^k)$ that is required to determine
$\nabla \varphi_\lambda(y^k)$ in line 3 of \Cref{alg:FB-LBFGS}.
The involved operations can be parallelized on a \ac{gpu} as we shall discuss in
Section \ref{sec:dual-gradient}, but the computations of $x(y)$, $x(d)$ and
$\nabla\hat{f}(y^k)R_\lambda(y^k)$ cannot be parallelized.

Often, the Lipshcitz constant of the gradient of $\hat{f}$ is
not known and needs to be estimated with a backtracking procedure.
The original backtracking proposed in~\cite{TheStePat16} halves the
value of $\lambda$ after the line search in line 4 if the following
condition is satisfied
\begin{multline}\label{eq:original-backtracking}
  \hat{f}(T_\lambda(w^k)) > \hat{f}(y^k) - \lambda \<\nabla \hat{f}(y^k), R_\lambda(y^k)\> \\
  + \tfrac{(1-\beta)\lambda}{2}\|R_\lambda(y^k)\|^2,
\end{multline}
for some $\beta\in[0, 1)$. The values $\hat{f}(y^k)$,
$\<\nabla \hat{f}(y^k), R_\lambda(y^k)\>$ and $\|R_\lambda(y)\|^2$
are known from the preceding line search, so the cost of the
backtracking is that of computing $\hat{f}(T_\lambda(w^k))$.
Alternatively, we may use the backtracking method proposed in
\cite[Linesearch 1]{BelloCruz2016} which halves $\lambda$ if
\begin{equation}\label{eq:simple-backtracking}
  \lambda\|\nabla \hat{f}(T_\lambda(y^k)) - \nabla \hat{f}(y^k)\|
  {}>{}
  \epsilon'' \|T_\lambda(y^k)-y^y\|,
\end{equation}
where $\epsilon''{}\in{}(0, \nicefrac{1}{2})$.
This backtracking procedure has a lower computational cost compared
to \Cref{eq:original-backtracking}.
In both cases, the L-BFGS buffer is emptied when the value of
$\lambda$ is updated.

\subsection{Parallelizable Newton-type Alternating Minimization Algorithm}%
\label{sec:parallelizable-nama}
The Newton-type alternating minimization algorithm (\nama) can be used to solve
the dual optimization problem $\mathbb{D}(p)$ in \Cref{eq:dual-problem} without
the need to compute the gradient of the FBE \cite{nama2019}.
\nama, applied to the dual optimization problem is given in \Cref{alg:nama}.

\nama involves a simple line search which consists in determining a $\tau_k$
so that the dual vector defined as $w^{k} = y^k + \tau_k d^k + (1-\tau_k)r^k$
satisfies the descent condition $\varphi_\lambda(w^{k})\leq \varphi_\lambda(y^k)$.
Again, if $\hat{f}$ is a quadratic function, we can precompute certain quantities
in a fashion akin to \Cref{eq:linesearch_trick}. In particular, before the
line search in line 7 of \Cref{alg:nama} we need to compute $x(r)$
and $x(d)$.

\begin{algorithm}[htb]
  \caption{\nama method for the dual optimization problem}\label{alg:nama}
  \begin{algorithmic}[1]
    \REQUIRE $\lambda\in(0,\mu_f/\|H\|^2)$,  $y^0$, $\epsilon>0$ (tolerance)
    \ENSURE Primal-dual solution triple $(x, z, y)$
    \STATE $k {}={} 0$
    \WHILE{$\|R_\lambda(y^k)\|>\epsilon$}
    \STATE $x^k {}={} x(y^k)$, $z^k {}={} z_\lambda(y^k)$
    \STATE $r^k {}={} z^k - Hx^k$
    \STATE $d^k {}={} -B^k r^k$ (Compute an L-BFGS direction using the L-BFGS buffer)
    \STATE Choose the smallest $\tau_k \in \{2^{-\nu}\}_{\nu\in\N}$ so that
    \begin{equation}\label{eq:decrease_condition_nama}
      \varphi_{\lambda}(w^{k})\leq \varphi_{\lambda}(y^k),
    \end{equation}
    where $w^{k} = y^k + \tau_k d^k + (1-\tau_k)r^k$
    \STATE $\tilde{x}^k {}={} x(w^k)$, $\tilde{z}^k {}={} z_\lambda(w^k)$
    \STATE $y^{k+1} {}={} y^k + \lambda(H \tilde{x}^k - \tilde{z}^k)$
    \STATE Compute the L-BFGS-related quantities $s^k = y^{k+1} - y^k$,
    $q^k = R_{\lambda}(y^{k+1}) - r^{k}$
    and $\rho_k = \langle s^k, q^k \rangle$
    \IF{$\<s^k, q^k\> > \epsilon' \|s^k\|^2 \|r^{k}\|^2$}
    \STATE Push $(s^k,q^k,\rho^k)$ into the L-BFGS buffer
    \ENDIF
    \STATE $k = k + 1$
    \ENDWHILE
    \STATE \textbf{return} $(x, z, y)=(x(y^k), z(y^k), y^k)$
  \end{algorithmic}
\end{algorithm}

The main computational cost involved in \Cref{alg:nama} comes from the
evaluation of $x(y)$, $x(r)$, and $x(d)$. Note that if $x$ is linear,
$x(w)$ can be computed at a very low computational cost.
In particular, the extrapolated vector $w^k$ can be written as
\(w^k = \tilde{y}^k + \tau_k \tilde{d}^k\), where \(\tilde{y}^k = y^k + r^k\)
and \(\tilde{d}^k = d^k - r^k\), therefore the decrease condition of
\nama in \Cref{eq:decrease_condition_nama} is equivalent to
\Cref{eq:linesearch_trick} with $\tilde{y}^k$ and $\tilde{d}^k$ in lieu
of $y^k$ and $d^k$ respectively, that is,
\begin{equation}\label{eq:linesearch_trick:2}
  \alpha_2(\tilde{y}, \tilde{d}) \tau^2
  +
  \alpha_1(\tilde{d}) \tau
  +
  \alpha_0(\tau; \tilde{y}, \tilde{d}) \leq 0.
\end{equation}

Overall, given that the computation of Hessian-vector products in \minfbe comes
at approximately the same cost as computing the dual gradient, and
given that the computation of $x(r)$ and $x(d)$ can be carried out
in parallel, \nama has a lower per-iteration computation cost.
Although \minfbe and \nama exhibit similar convergence properties,
with \nama we can afford a greater parallelizability that leads to
superior performance in practice as we shall show in \Cref{sec:numerical_simulations}.

\subsection{Efficient parallel computations}\label{sec:dual-gradient}
\acp{gpu} have a hardware architecture that allows the execution of the same
set of instructions on different memory positions. \acp{gpu} are equipped
with a set of SIMD stream processors, each having its own
computing resources, that execute ``compute kernels,'' that is, functions
that are executed simultaneously on different data.

NVIDIA's \acp{gpu} use the CUDA programming interface
where kernels are executed in parallel threads, which
are organised in \textit{blocks} which can share memory and which are
in turn organised in grids. At a hardware level, threads are
executed in parallel in \textit{warps} of 32 threads.
Threads in the same block have asynchronous
read/write access to a local shared memory and can synchronize.
Each thread has its own local memory, and all threads have access to
the device's global memory.
Modern \acp{gpu} count several streaming multiprocessors with hundreds of cores,
possess a computing throughput of several Tera-FLOPs, and have a
significant memory capacity of several GBs.
The hardware architecture and programming model of \acp{gpu} necessitates
a fresh look at parallelization approaches for numerical optimization.
Kernels are best suited for the parallel execution of simple numerical
operations.

The efficient computation of the dual gradient is of crucial
importance for the performance of the algorithm we are about to
describe. By virtue of the Conjugate Subgradient
Theorem~\cite[Theorem 23.5]{rockafellar-convex-analysis},
we have that
\begin{multline}\label{eq:nabla-f-star}
  x(y)
  {}={}
  \argmin_{z\in\mathcal{Z}(p)} \bigg\{
  \sum_{i\in\nodes(1, N)}\hspace{-1.2em}
  \pi^i\hat\phi^i(x^{\anc(i)},u^{\anc(i)})\notag
  \\[-0.5em]
  {}+{}
  \sum_{i\in\nodes(N)}\hspace{-0.9em}\pi^i \hat\phi_N(x^i)\bigg\},
\end{multline}
where $\hat\phi^i(x^{\anc(i)},u^{\anc(i)}) = \phi^i(x^{\anc(i)},u^{\anc(i)})
  + \<y^i, F^i x^{\anc(i)}+ G^{i}u^{\anc(i)}\>$ for $i\in\nodes(1, N)$, and
$\hat\phi_N(x^i)=\phi_N(x^i)+\<y^i, F_N^i x^i\>$, for $i\in\nodes(1, N)$.
The solution of this problem can be determined via a dynamic programming
in a way akin to~\cite[Algorithm~1]{SamSopBemPat15} leading to
\Cref{alg:solve} wherein $\Phi_k^i$, $\Theta_k^i$, $D_k^i$, $\Lambda_k^i$, $K_k^i$
$\sigma_k^i$, $c_k^i$ are computed once offline following a Riccati-type recursion.
In cases where the data of the optimal control problem need to be updated (e.g.,
if the dynamical system is time varying, or the parameters of the cost must be
updated in real time), the computation of these matrices can be carried out on
a \ac{gpu} and in fact the time for their computation is negligible compared to
that of solving the problem.

\begin{algorithm}[htbp]
  \caption{Computation of the dual gradient, $x(y)$}
  \label{alg:solve}
  \begin{algorithmic}[1]
    \REQUIRE Dual vector $y\in\R^m$
    \ENSURE $x(y)$
    \STATE $\hat{q}^i\gets y^i,$ for all $i\in\nodes(N)$
    \FOR[\texttt{in parallel}]{$i\in\nodes(0, N-1)$}
    \STATE  $\phantom{\hat{q}^i}\mathllap{u^i}\gets\Phi^i y^i +\sigma^i$
    \STATE $\hat{q}^i\gets  D^{i\top} y^i + \hat{c}^i$
    \ENDFOR
    \FOR{$k = N - 1, \ldots, 0$}
    \FOR[\texttt{in parallel}]{$i\in\nodes(k)$}
    \STATE  $\phantom{\hat{q}^i}\mathllap{u^i}\gets \sum_{i_\smallplus \in \child(i)}\Theta^{i_\smallplus} \hat{q}^{i_\smallplus}$
    \STATE $\hat{q}^i\gets  \sum_{i_\smallplus \in \child(i)} \Lambda^{i_\smallplus\top} \hat{q}^{i_\smallplus}$
    \ENDFOR
    \ENDFOR
    \STATE $x^0=p$
    \FOR{$i\in\nodes(0, N-1)$}
    \STATE $u^i\gets K^ix^i+u^i$
    \FOR{$i_\smallplus \in \child(i)$} 
    \STATE $x^{i_\smallplus}\gets A^{i_\smallplus} x^i+B^{i_\smallplus} u^i +c^{i_\smallplus}$
    \ENDFOR
    \ENDFOR
    \RETURN $x(y) = (\{x^i\},\{u^i\})$.
  \end{algorithmic}
\end{algorithm}

The computation of Hessian-vector products of the form
$\nabla^2 \hat{f}(y){}\cdot{}r$ that is required for the
computation of the gradient of the FBE is given in
\Cref{alg:hessian-vec-products}. \Cref{alg:solve} and
\Cref{alg:hessian-vec-products} incur roughly the same
computation cost.
\begin{algorithm}[ht]
  \caption{Computation of Hessian-vector products required for the
    computation of $\nabla \varphi_\lambda$}\label{alg:hessian-vec-products}
  \begin{algorithmic}[1]
    \REQUIRE Vector $r$
    \ENSURE Hessian-vector product, $\nabla^2 \hat{f}(y){}\cdot{}r$
    \STATE $\hat{q}^i\gets r^i,$ for all $i\in\nodes(N)$
    \FOR{$k = N - 1, \ldots, 0$}
    \FOR[\texttt{in parallel}]{$i\in\nodes(k)$}
    \STATE  $\phantom{\hat{q}^i}\mathllap{\hat{u}^i}\gets \Phi^i r^i + \sum_{i_\smallplus \in \child(i)}\Theta^{i_\smallplus} \hat{q}^{i_\smallplus}$
    \STATE $\hat{q}^i\gets  D^{i\top} r^i + \sum_{i_\smallplus \in \child(i)} \Lambda^{i_\smallplus\top} \hat{q}^{i_\smallplus}$
    \ENDFOR
    \ENDFOR
    \STATE $\hat{x}^0=0$
    \FOR{$i\in\nodes(0, N-1)$}
    \STATE $\hat{u}^i\gets K^i\hat{x}^i+\hat{u}^i$
    \FOR{$i_\smallplus \in \child(i)$} 
    \STATE $\hat{x}^{i_\smallplus} {}\gets{} A^{i_\smallplus} \hat{x}^i+B^{i_\smallplus} \hat{u}^i$
    \ENDFOR
    \ENDFOR
    \RETURN $\nabla^2 \hat{f}(y){}\cdot{}r {}={} (\{\hat{x}^i\},\{\hat{u}^i\})$.
  \end{algorithmic}
\end{algorithm}

Lastly, most proximal operations can be massively parallelized.
For example, if $\bar\phi^{i}(z) = \delta(z {}\mid{} Y^i)$ and
$Y^i=\{z{}:{}z^{i}_{\min} \leq z \leq z^i_{\max}\}$, then the
computation of $\prox_{\lambda \bar{\phi}^i} = \proj_{Y^i}$
is element-wise independent and can be easily parallelized.
Likewise, a great many proximal operators, such as those of
the indicators of rectangles and common norm-balls, and functions
such as $\|{}\cdot{}\|_1$, the Huber loss function and more,
lend themselves to high parallelizability \cite{parikh-boyd-proximal}.

In general, the total memory that needs to be allocated on the \ac{gpu}
grows linearly with the length of the L-BFGS buffer, linearly with the
prediction horizon, and linearly with the number of nodes of the tree,
and quadratically with the system states and inputs.
The additional parallelisation in \nama requires the allocation of
additional memory on the \ac{gpu}, but leads to a higher throughput and
occupancy of the device.

\subsection{Preconditioning}
Stochastic optimal control problems tend to be ill conditioned
because of the presence of generally small probability values.
As first-order methods are known to be affected by the problem being
ill conditioned, here we make use of a simple diagonal preconditioning heuristic
where we scale the original dual variables $y=((y^i)_{i\in\nodes(0, N)},
  (y_{N}^{i})_{i\in\nodes(N)})$ by introducing the scaled
dual variables $\bar{y}=((\bar{y}^i)_{i\in\nodes(0, N)}, (\bar{y}_N^i)_{i\in\nodes(N)})$ with
\begin{equation}
  \bar{y}^i = \frac{y^i}{\sqrt{\pi^i}},
\end{equation}
for $i\in\nodes(0, N)$ and
\begin{equation}
  \bar{y}_N^i = \frac{y_N^i}{\sqrt{\pi^i}},
\end{equation}
for $i\in\nodes(N)$.
This scaling is a heuristic similar to
the Jacobi preconditioning discussed in \cite{Giselson2015}.

\subsection{Warm start}
Generally, the accelerated projected gradient method converges at a rate $\mathcal{O}(1/k^2)$
and although it may exhibit slow convergence, its iterations are computationally cheap,
so it can be used to warm start \minfbe and \nama.
We have observed that running as few as five iterations of
\gpad \cite{PatrinosBemporad2014,SamSopBem+2018} can provide a
good warm starting point for \minfbe and \nama.


\section{Numerical Simulations}\label{sec:numerical_simulations}
This section is organised in two parts: in \Cref{sec:spring-mass}
we compare \minfbe and \nama with the accelerated projected gradient
method and discuss the convergence rate of each method.
In particular, we demonstrate that a serial implementation of \nama
and \minfbe leads to superior performance compared to the accelerated
proximal gradient method. The two methods exhibit comparable convergence
speed.
Next, in \Cref{sec:barcelona-dwn} we apply \minfbe and \nama to solve
a large-scale stochastic optimal control problem for the operating
management of the drinking water network of Barcelona taken from \cite{SamSopBem+2018}.
We show that \nama affords a higher parallelisation leading to a
significant performance improvement.

\subsection{Spring-mass-damper array}\label{sec:spring-mass}
\begin{figure}[htbp]
   \centering
   \includegraphics[width=0.98\linewidth]{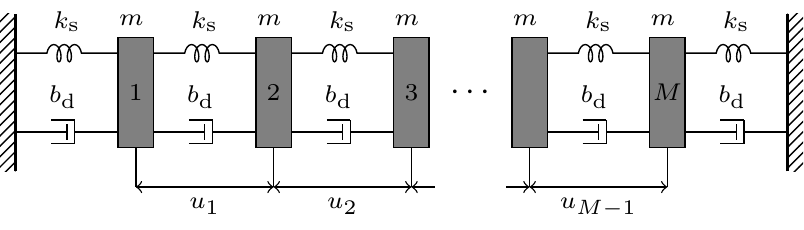}
   \caption{Array of $M$ consecutive interconnected masses, $M+1$ elastic springs and dampers
      and $M-1$ actuators.}
   \label{fig:spring_mass}
\end{figure}

Consider an array of $M$ consecutive point particles of mass $m$ connected to each other through
elastic springs of stiffness $k_{\rm s}$ and linear dampers with viscous damping
coefficients $b_{\rm d}$ illustrated in \Cref{fig:spring_mass}.

In between the successive masses $j$ and $j+1$, for $j=1,\ldots,M-1$,
there is an actuator that can apply a force $u_j\in[u_{\rm min}, u_{\rm max}]$.
The state variable of this sytem comprises of the positions $p_j$ of the masses
and their velocities $v_j$, which are constrained in $[p_{\rm min}, p_{\rm max}]$
and $[v_{\rm min}, v_{\rm max}]$, respectively.
The system is described by a set of linear differential equations which can be obtained
by the application of Newton's second law of motion, which, after discretisation
with sampling time $T_s$ and a zero-order hold, yields a discrete-time linear
time invariant system.
Furthermore, we assume that there is an external additive disturbance $c_{w_k}$, as in
\Cref{eq:affine_stochastic_dynamics}, which is driven by a discrete Markov process, $w_k$,
with two modes.

In this example, we consider a stochastic optimal control problem with
prediction horizon $N$, quadratic stage cost functions
$\phi^i(x, u) = x^\top Q x + u^\top R u$,
and quadratic terminal costs
$\phi_N^i(x)=x^\top Q_N x$.
Moreover, we have $M=5$ masses with
\(m=\unit[5]{kg}\),
\(k_{\rm s}=\unitfrac[1]{N}{m}\),
\(b_{\rm d}=\unitfrac[0.1]{Ns}{m}\),
\(u_{\rm max} = -u_{\rm min} = \unit[2]{N}\)
and the maximum allowed velocity is
\(\unitfrac[5]{m}{s}\).
The prediction horizon is $N=11$ and the external disturbance $c_k$ is driven by a Markov
chain with two modes with initial probability distribution $p_c=(0.5, 0.5)$
and probability transition matrix $P_c=\smallmat{0.1 & 0.9\\0.9 & 0.1}$;
at mode 1 the value of $c$ is zero and at mode 2, $c$ takes the value $0.1$.
The sampling time is $T_s=\unit[0.5]{s}$. Lastly, the weights of the stage and terminal
cost functions are $Q=5I_{10}$, $R=2I_{4}$ and $Q_N = 100I_{10}$.
No warm starting is used in any of the algorithms.

\begin{figure}[!htbp]
   \centering
   \includegraphics[width=0.96\linewidth]{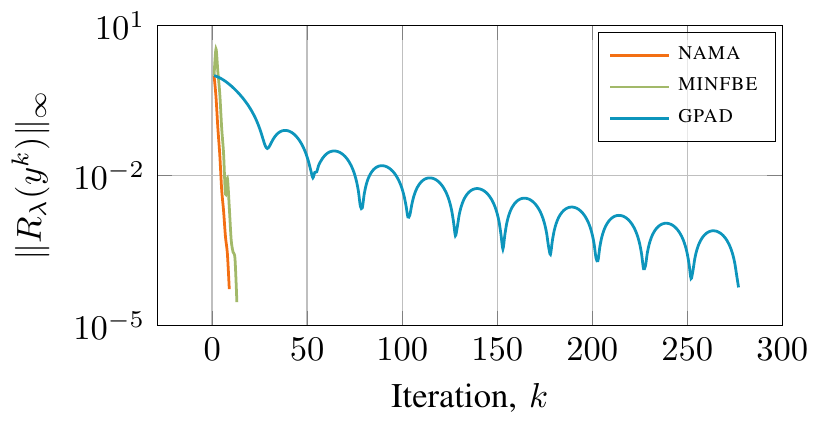}
   \caption{Comparison of the convergence of \nama, \minfbe, and the accelerated projected gradient (\gpad) method
      applied to the dual problem.}
   \label{fig:Convergence_sample_point}
\end{figure}

\begin{figure}[t]
   \centering
   \includegraphics[width=0.98\linewidth]{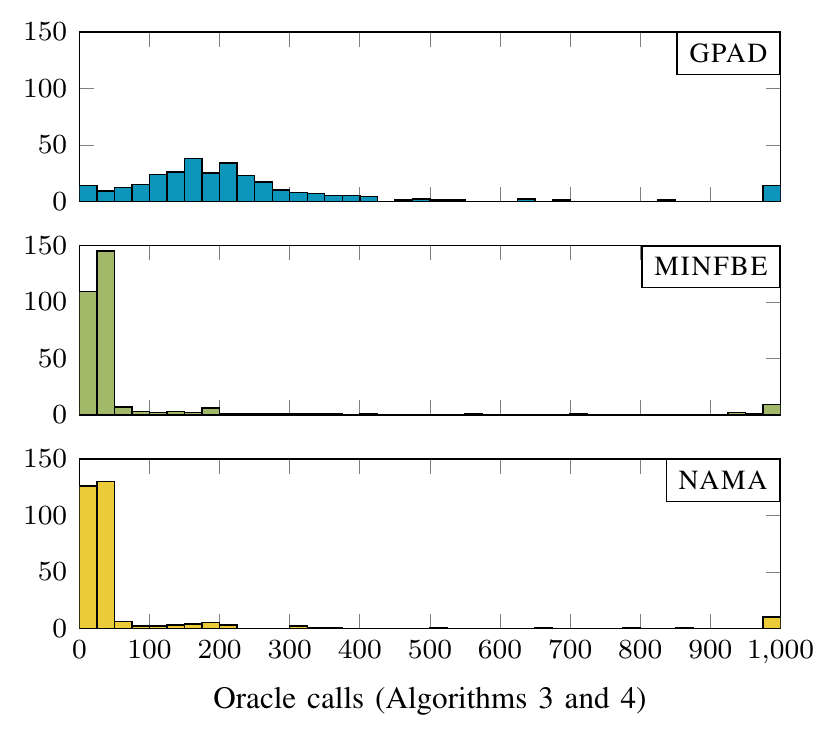}
   \caption{Distribution of the number of oracle calls required for the computation of the dual gradient (\Cref{alg:solve})
      and Hessian-vector products (\Cref{alg:hessian-vec-products}) for \textsc{gpad}, \minfbe and \nama.}
   \label{fig:histogram_hessian_call}
\end{figure}

We ran the stochastic optimal control problem for $300$ initial states $x_0=p$, sampled uniformly
from the problem's domain. These problems were solved with \nama, \minfbe and the accelerated projected
gradient method applied to the dual problem (\gpad) following \cite{SamSopBem+2018}.
In \nama and \minfbe we used L-BFGS directions with a memory of $5$.
We used the same termination condition in all methods with $\epsilon=5\cdot 10^{-4}$.

\gpad is known to converge at a rate of $\mathcal{O}(1/k^2)$, which can be observed in \Cref{fig:Convergence_sample_point};
clearly, \gpad can only achieve low to medium accuracy solutions within a few hundred iterations.
On the other hand, \minfbe and \nama exhibit a significantly faster convergence rate and require fewer
iterations to achieve solutions of higher accuracy.

In \Cref{fig:histogram_hessian_call} we show the number of calls of \Cref{alg:solve} and
\Cref{alg:hessian-vec-products} required to solve the aforementioned collection of 300 random
problems up to the desired accuracy. We may observe that in the majority of cases ($84\%$),
\minfbe and \nama can solve the problems with no more than 50 calls, whereas the median
of the number of calls corresponding to \gpad is 188.

Note that \nama and \minfbe appear to perform on a par. However, in the next section
we will demonstrate that \nama allows for greater parallelizability leading to superior
performance on a \ac{gpu}.

\subsection{Large-scale drinking water network}\label{sec:barcelona-dwn}
In this section we apply the proposed numerical optimization methods for the
solution of a model predictive control problem for a drinking water network,
whose transportation dynamics is described by
\begin{subequations}\label{eq:dwn:mass_balance}
   \begin{align}
      x_{t+1} {}={} & Ax_t + Bu_t + G d_{t},
      \label{eq:dwn:mass_balance::1}
      \\
      0 {}={}       & E_u u_t + E_d d_{t},
      \label{eq:dwn:mass_balance::2}
   \end{align}
\end{subequations}
where $x_t$ is the vector of the volume of water in the reservoirs of the
network, $u_t$ is the vector of pumping set points and $d_t$ is the vector
of water demands from the various distribution nodes.
The value of $d_t$ is measured at time $t$ and future demand values are
predicted by a model that returns estimates $\hat{d}_{t+t'{}\mid{}t}$,
for $t' \geq t$, while $d_{t+t'} = \hat{d}_{t+t'{}\mid{}t} + \epsilon_{t'}$,
where $\epsilon_{t'}$ is a random process that can be described by a
scenario tree~\cite{SamSopBem+2018}.

The water network model \eqref{eq:dwn:mass_balance} comprises
$63$ states corresponding to water level in the tanks, $114$ inputs
corresponding to flow control devices (pumps and valves),
$88$ disturbance variables corresponding to the demand sectors and
input-disturbance relationship corresponding to the 17 mixing nodes.
The detailed stochastic optimal control problem and the formulation
of the optimisation problem is discussed in~\cite{SamSopBem+2018}.
The operation of the water network is subject to uncertainty in
water demand and electricity prices.

\begin{figure}
   \centering
   \includegraphics[width=0.96\linewidth]{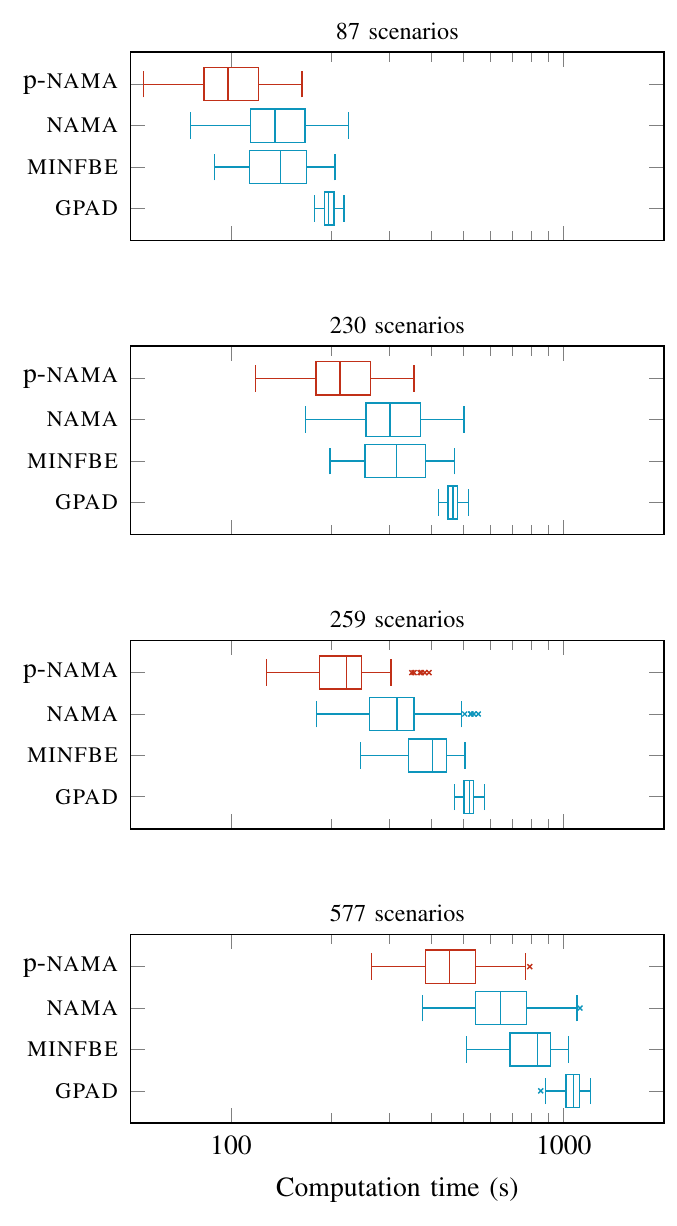}
   \caption{Box plots of the computational time with different scenario tree sizes
      for the \ac{gpu} implementation of the algorithms:
      (parallel) \nama, \minfbe and \gpad.}
   \label{fig:Computation_time}
\end{figure}

The \nama and \minfbe algorithms are implemented in the
\texttt{RapidNet}\footnote{\url{https://github.com/GPUEngineering/RapidNet}} software
package that is developed for the operational control of water network problems.
All simulations presented in this section were carried out on an NVIDIA Tesla C2075
\ac{gpu} which counts 448 CUDA cores running at $\unit[1.15]{GHz}$ and $\unit[6]{GB}$
of dedicated memory.

In order to demonstrate the effect of the additional parallelisation
in \nama that we discussed in \Cref{sec:parallelizable-nama},
we provide results for the method with that additional parallelization
in the computation of the line search (p-\nama) and \nama without that additional parallelization.

The parallel computations involved in Algorithms \ref{alg:solve} and \ref{alg:hessian-vec-products}
are carried out using cuBLAS's \texttt{cublasSgemmBatched} and \texttt{cublasSgemm}.
In this example, matrices $A$, $B$, $G_d$, $E_u$ and $E_d$ are sparse,
and this has been used to tailor the implementations of Algorithms
\ref{alg:solve} and \ref{alg:hessian-vec-products} to be more
efficient.


The L-BFGS memory is set to $15$. In the case with 577 scenarios,
the problem involves 2.1 million primal and 3.8 million dual variables and
\nama and \minfbe algorithms require an excess of
$2.9\%$  ($\unit[172]{MB}$) of memory and p-\nama requires an excess
of $4.1\%$ ($\unit[242]{MB}$) of memory than dual accelerated proximal gradient (\gpad) algorithm.
The solve times of p-\nama, \nama, \minfbe and \gpad are shown in
\Cref{fig:Computation_time}, where note that the horizontal axis is logarithmic.
It can be observed that p-\nama is noticeably faster compared to \nama, \minfbe and
\gpad.



\section{Conclusions}
In this paper we proposed the use of \minfbe and \nama for solving large-scale
scenario-based convex stochastic optimal control problems. Both methods use
limited-memory quasi-Newtonian, L-BFGS, directions and exhibit a very fast
convergence rate. They are both suitable for parallelization on GPUs, but
\nama lends itself to a significantly higher parallelization.
We presented compelling results on two stochastic optimal control
problems, namely a spring-mass-damper array and the drinking water network
of Barcelona, demonstrating that the two methods significantly outperform \gpad,
whose parallelizable implementation on a \ac{gpu} has been previously shown to
outperform Gurobi's interior point solver \cite{SamSopBem+2018}.
Future work will focus on the development of parallelizable methods for
large-scale scenario-based risk-averse optimal control problems
\cite{risk-averse-mpc-2019}.

\bibliographystyle{IEEEtran}
\bibliography{bibliography}

\end{document}